\documentstyle[12pt]{article}
\begin{document}
\begin{center}
\today\\[10pt]
{\Large\bf Some Applications of the Ultrapower Theorem to the Theory of 
Compacta}
\\[20pt]
Paul Bankston\\
Department of Mathematics, Statistics and Computer Science\\ 
Marquette University\\
Milwaukee, WI 53201-1881\\[20pt]

{\it Dedicated to Bernhard Banaschewski on the occasion of his 70th
Birthday.}
\end{center}
\begin{abstract}
The ultrapower theorem of Keisler-Shelah allows such model-theoretic notions
as elementary equivalence, elementary embedding and existential embedding to be
couched in the language of categories (limits, morphism diagrams).  This in 
turn allows analogs of these (and related) notions to be transported into
unusual settings, chiefly those of Banach spaces and of compacta.
Our interest here is the enrichment of the theory of compacta, especially
the theory of continua, brought about by the immigration of model-theoretic
ideas and techniques.
  
\end{abstract}

\noindent
{\it A.M.S. Subject Classification\/} (1991):
03C20, 54B35, 54C10, 54D05, 54D30, 54D80, 54F15, 54F45, 54F55.\\[10pt]


\section{Introduction.}\label{1}

In the ten years since the celebration of Professor Banaschewski's
60th birthday (and the writing of \cite{Ban4} in commemoration), there
has been a fair amount of development in the theory of ultracoproducts
of compacta (i.e., compact Hausdorff spaces).
Papers \cite{Ban3, Ban5, Ban7, Ban8} have been written
by this author; also there is the paper \cite{Gur} by the late R.
Gurevi\v{c}.  In addition to this, there has been a parallel development
in the first-order theory of Banach spaces/algebras begun by C. W. Henson
in \cite{Hen1}.  There are important links between Banach space theory
and our work (mainly through Gel'fand-Na\u{\i}mark duality), and we 
present here two applications of Banach techniques to the theory of compacta. 
(See \ref{3.1} and \ref{4.2} below.)
Since Gel'fand-Na\u{\i}mark duality is a two-way street, much of the
``dualized'' model theory developed for compacta may be directly translated
into the Banach model theory of commutative $B^*$-algebras. (See
the results in \S5 and \S6 below.)  It seems likely that the future will 
see much in the way of progress in these two streams of research, as
topological issues stimulate the analytic and {\it vice versa}.

We begin with a quick review of the topological ultracoproduct construction;
detailed accounts may be found in \cite{Ban2, Ban3, Ban4, Ban5, Ban7, Ban8,
Gur}. 

We let {\bf CH} denote the category of compacta
and continuous maps.  In model theory, it is well known
that ultraproducts (and reduced products in general, but we restrict
ourselves to maximal filters on the index set) 
may be described in the language of category theory; i.e., as direct limits 
of (cartesian) products, where the directed 
set is the ultrafilter with reverse inclusion, and the system of products  
consists of cartesian products taken over the various sets in the ultrafilter.
(Bonding maps are just the obvious restriction maps.)  When we transport this
framework to {\bf CH}, the result is somewhat less than spectacular:  If
$\langle X_i:i \in I\rangle$ is a family of nonempty compacta and $\cal D$
is a nonprincipal ultrafilter on $I$, then the {\bf CH}-ultraproduct is 
degenerate (i.e., has only one point).  
What turns out to be vastly more fruitful, however, is the ultraproduct 
construction in the category-opposite of {\bf CH}; i.e., take an inverse limit 
of coproducts.  The result is the {\bf topological ultracoproduct}, 
and may be concretely
(if opaquely) described as follows:  Given $\langle X_i:i  \in I \rangle$ and
$\cal D$, let $Y$ be the disjoint union $\bigcup_{i \in I}(X_i \times \{i\})$
(a locally compact space).  With $q:Y \to I$ the natural projection onto
the second co\"{o}rdinate (where $I$ has the discrete topology), we then
have the Stone-\v{C}ech lifting $q^{\beta}:\beta(Y)\to \beta(I)$.  Now the
ultrafilter $\cal D$ may be naturally viewed as an element of $\beta(I)$, 
and it is not hard to show that the topological ultracoproduct 
$\sum_{\cal D}X_i$ is the pre-image 
$(q^{\beta})^{-1}[\cal D]$.  
(The reader may be familiar with the Banach ultraproduct \cite{DK}. 
This construction is indeed the ultraproduct
in the category of Banach spaces and nonexpansive linear maps, and may
be telegraphically described using the recipe: take the usual ultraproduct,
throw away the infinite elements, and mod out by the subspace of 
infinitesimals.  Letting $C(X)$ denote the Banach space of continuous
real-valued (or complex-valued) continuous functions with $X$ as domain,
the Banach ultraproduct of $\langle C(X_i): i\in I\rangle$ via $\cal D$ is
just $C(\sum_{\cal D}X_i)$.)

If $X_i = X$
for all $i\in I$, then we have the {\bf topological ultracopower}
$XI\backslash {\cal D}$, a
subspace of $\beta(X\times I)$. In this case there is the Stone-\v{C}ech
lifting $p^{\beta}$ of the natural first-co\"{o}rdinate map 
$p:X\times I\to X$.  Its restriction to the ultracopower is a continuous
surjection, called the {\bf codiagonal map}, and is officially 
denoted $p_{X,{\cal D}}$ (with the occasional notation-shortening alias
possible). 
This map is dual to the natural diagonal map from a relational structure
to an ultrapower of that structure, and is 
not unlike the standard part map from nonstandard
analysis.)

Many notions from classical first-order model theory, principally elementary
equivalence, elementary embedding and existential embedding, may be phrased in 
terms of mapping
conditions involving the ultraproduct construction.  Because of the 
(Keisler-Shelah) ultrapower theorem (see, e.g., \cite{CK}), two relational
structures are elementarily equivalent if and only if some ultrapower of one
is isomorphic to some ultrapower of the other; a function from one relational
structure to another is an elementary embedding if and only if there is an
ultrapower isomorphism so that the obvious square mapping diagram commutes.
A function $f:A \to B$ between relational structures is an existential embedding
(i.e., making the image under $f$ a substructure of $B$ that is existentially
closed in $B$) if and only if there are embeddings $g:A \to C$, $h:B \to C$
such that $g$ is elementary and equal to the composition $hf$.
($C$ may be taken to be an ultrapower of $A$, with $g$ the natural diagonal.)
(see also, e.g., \cite{Ban2,Ban5,Ekl}).  

In {\bf CH} one then constructs ultracoproducts,
and talks of co-elementary equivalence, co-elementary maps and 
co-existential maps. 
Co-elementary
equivalence is known \cite{Ban2,Ban5,Gur} to preserve important properties
of topological spaces, such as being infinite, being a continuum (i.e., 
connected), being Boolean (i.e., totally
disconnected), having (Lebesgue) covering dimension $n$, and being a
decomposable continuum.  If $f:X \to Y$ is a co-elementary map in {\bf CH},
then of course $X$ and $Y$ are co-elementarily equivalent ($X \equiv Y$). 
Moreover, since $f$ is a continuous surjection (see 
\cite{Ban2}), additional information about $X$ is transferred to $Y$.  For
instance, continuous surjections in {\bf CH} cannot raise {\bf weight\/} (i.e., 
the smallest cardinality of a possible topological base, and for many
reasons the right cardinal invariant to replace cardinality in the dualized
model-theoretic setting), so metrizability
(i.e., being of countable weight in the compact Hausdorff context) is
preserved.  Also local connectedness is preserved, since continuous surjections
in {\bf CH} are quotient maps.  Neither of these properties is an invariant
of co-elementary equivalence alone.  

When attention is restricted to the full subcategory of Boolean 
spaces, the dualized model theory matches perfectly with
the model theory of Boolean algebras because of Stone duality.  In the
larger category there is no such match \cite{Bana,Ros}, however, and one is
forced to look for other (less direct) model-theoretic aids.  Fortunately  
there is a finitely axiomatizable AE Horn class of bounded distributive 
lattices,
the so-called {\bf normal disjunctive\/} lattices \cite{Ban8} (also called
Wallman lattices in \cite{Ban5}), comprising precisely the (isomorphic
copies of) {\bf lattice bases}, those lattices that serve as bases for the 
closed sets of compacta.  (To be more specific:  The normal disjunctive
lattices are precisely those bounded lattices $A$ such that there exists
a compactum $X$ and a meet-dense sublattice ${\cal A}$ of the closed set
lattice $F(X)$ of $X$ such that $A$ is isomorphic to $\cal A$.)  
We go from lattices to spaces, as in the case of Stone
duality, via the {\bf maximal spectrum\/} $S(\;)$, pioneered by H. Wallman 
\cite{Walm}.
$S(A)$ is the space of 
maximal proper filters of $A$; a typical basic closed set in $S(A)$ is the
set $a^{\sharp}$ of elements of $S(A)$ containing a given element $a \in A$.  
$S(A)$ is generally compact with this topology.
Normality, the condition that if
$a$ and $b$ are disjoint ($a \sqcap b = \bot$), then there are $a'$,
$b'$ such that $a \sqcap a' = b \sqcap b' = \bot$ and $a' \sqcup b' = \top$,
ensures that the maximal spectrum topology is Hausdorff.  Disjunctivity,
which says that for any two distinct lattice elements there is a nonbottom 
element that is below one of the first two elements
and disjoint from the other,
ensures that the map $a \mapsto a^{\sharp}$ takes
$A$ isomorphically onto the canonical closed set base for $S(A)$.  $S(\;)$
is contravariantly functorial: If $f:A \to B$ is a homomorphism of normal
disjunctive lattices and $M \in S(B)$, then $f^S(M)$ is the unique maximal
filter extending the prime filter $f^{-1}[M]$.  (For normal 
lattices, each prime filter is contained in a unique maximal one.)
It is a fairly
straightforward task to show, then, that $S(\;)$ converts ultraproducts
to ultracoproducts, elementarily equivalent lattices to co-elementarily
equivalent compacta, and elementary (resp., existential) embeddings to 
co-elementary (resp., co-existential)
maps.  Furthermore, if $f:A\to B$ is a {\bf separative} embedding; i.e.,
an embedding such that
if $b \sqcap c = \bot$ in $B$, then there exists $a \in A$ such that
$f(a) \geq b$ and $f(a) \sqcap c = \bot$, then $f^S$ is a homeomorphism  
(see \cite{Ban2,Ban4,Ban5,Ban8,Gur}).  Because of this, there is much
flexibility in how we may obtain $\sum_{\cal D}X_i$:  Simply choose a
lattice base ${\cal A}_i$ for each $X_i$ and apply $S(\;)$ to the
ultraproduct $\prod_{\cal D}{\cal A}_i$.

\section{The Topological Behavior of Co-existential Maps.}\label{2}

Recall that a first-order formula in prenex form is an 
{\bf existential formula} if all its quantifiers are existential.
A function $f:A \to B$ between structures is an 
{\bf elementary (resp., existential) embedding} if for every formula 
(resp. existential formula)
$\varphi(x_1,...,x_n)$ and every
$n$-tuple $\langle a_1,...,a_n\rangle$ from $A$, $A \models 
\varphi[a_1,...,a_n]$ if and only if $B \models \varphi[f(a_1),...,f(a_n)]$.
 
The ultrapower theorem states that a function $f:A\to B$ is an elementary
embedding if and only if there is an isomorphism of ultrapowers
$h:A^I/{\cal D} \to B^J/{\cal E}$ such that the obvious mapping square
commutes; i.e., such that $d_{\cal E}f = hd_{\cal D}$, where $d_{\cal D}$
and $d_{\cal E}$ are the natural diagonal embeddings.  
There is also a characterization of existential embeddings along similar lines:
$f:A\to B$ is an existential embedding if and only if there are embeddings 
$g:A \to C$ and $h:B \to C$ such that $g$ is elementary and $g = hf$.
(By the ultrapower theorem, we may take $C$ to be an ultrapower $A^I/{\cal D}$
and $g = d_{\cal D}$.)     

These characterizations have inspired the {\it definition} of the
notions of co-elementary map and co-existential map in {\bf CH} in terms of
analogous mapping diagrams involving ultracopowers.  Co-elementary
(indeed, co-existential) maps are clearly continuous surjections.
However, since
specific situations involving co-elementary and co-existential maps are not 
guaranteed
to correspond to analogous situations in elementary classes of relational
structures, one may not take too much for granted.  For example, there is no
assurance {\it a priori} that the classes of co-elementary and co-existential
maps are closed under composition or terminal factors (i.e., if both $f$ 
and $gf$ have a particular property, then $g$ has the property).  As it
happens, it is shown in \cite{Ban2} that closure under composition and
terminal factors holds for co-elementary maps; also in \cite{Ban7} it is 
shown that there is an amalgamation property for co-elementary maps:  If 
$f:Y\to X$ and $g:Z\to X$ are co-elementary, then there exist co-elementary
maps $u:W \to Y$ and $v:W \to Z$ such that $fu = gv$.  We do not know whether
there is amalgamation for co-existential maps; however we can show closure 
under composition and terminal factors.  We first prove a useful lemma.\\

\subsection{Lemma.}\label{2.0.5} Let $\{f_{\delta}:X_{\delta} \to 
Y_{\delta}: \delta \in \Delta\}$ be a
family of co-elementary and co-existential maps between compacta.  
Then there is a
single ultrafilter witness to the fact. More precisely, there is an ultrafilter 
$\cal D$ on a
set $I$ such that: $(i)$ if $f_{\delta}$ is co-elementary, then there is a
homeomorphism $h_{\delta}: X_{\delta}I\backslash {\cal D} \to 
Y_{\delta}I\backslash {\cal D}$ 
such that $f_{\delta}p_{X_{\delta},{\cal D}} = 
p_{Y_{\delta},{\cal D}}h_{\delta}$; and $(ii)$ if
$f_{\delta}$ is co-existential, then there is a continuous surjection 
$g_{\delta}:Y_{\delta}I\backslash {\cal D} \to X_{\delta}$ such that 
$f_{\delta}g_{\delta} = p_{Y_{\delta},{\cal D}}$.\\

\noindent
{\bf Proof.}  In S. Shelah's (GCH-free)
proof of the ultrapower theorem (see \cite{She}), he proves that if 
$\lambda$ is an infinite cardinal and $\mu = \mbox{min}\{\kappa:
\lambda^{\kappa} > \lambda\}$, then there is an ultrafilter $\cal D$ on
$\lambda$ such that whenever $A$ and $B$ are elementarily equivalent
relational structures of cardinality $< \mu$, then $A^{\lambda}/{\cal D}$
and $B^{\lambda}/{\cal D}$ are isomorphic.  (We may easily extend this
theorem to cover the situation involving elementary embeddings by adding
constants to name domain elements.  This way we get commuting square
diagrams.)

Assume the maps $f_{\delta}$ are all co-elementary.  (Co-existential maps are
handled in a similar way.)  Then we have a collection of homeomorphisms
$k_{\delta}:X_{\delta}I_{\delta}\backslash {\cal D}_{\delta} \to 
Y_{\delta}J_{\delta}\backslash {\cal E}_{\delta}$ making
the obvious mapping squares commute.  Now, using Shelah's theorem, 
find an ultrafilter $\cal D$, on a very large index set $I$, so that the
elementarity of the diagonal embeddings $d_{\delta}:F(X_{\delta}) \to 
F(X_{\delta})^{I_{\delta}}/{\cal D}_{\delta}$
and $e_{\delta}:F(Y_{\delta}) \to 
F(Y_{\delta})^{J_{\delta}}/{\cal E}_{\delta}$ (between normal disjunctive
lattices), $\delta \in \Delta$, is simultaneously witnessed.  When we apply the
spectrum functor $S(\;)$ to these squares, we then have witnesses to
the co-elementarity of the co-diagonal maps $p_{\delta}: 
X_{\delta}I_{\delta}\backslash {\cal D}_{\delta}\to X_{\delta}$
and $q_{\delta}: Y_{\delta}I_{\delta}\backslash {\cal E}_{\delta}
\to Y_{\delta}$.  (Recall that $S(\;)$ converts canonical ultrapower diagrams
to canonical ultracopower diagrams, and that there is an iteration theorem
for ultracopowers: ``ultracopowers of ultracopowers are ultracopowers,''
see \cite{Ban2}.)
Now $(\;)I\backslash {\cal D}$, as an operator on compacta, 
is functorial;
it is covariant, and preserves (reflects as well) 
many interesting properties of continuous maps.    
In particular, as is straightforward to check, it preserves and reflects the 
properties of being surjective, of being one-one, of being co-elementary,
and of being co-existential.  Thus, when we apply $(\;)I\backslash {\cal D}$
to the homeomorphisms $k_{\delta}$, we have a new collection of homeomorphisms;
hence we have a homeomorphism linkage between 
$X_{\delta}I\backslash {\cal D}$ and
$Y_{\delta}I\backslash {\cal D}$, witnessing the co-elementarity of 
$f_{\delta}$.
$\dashv$\\

\subsection{Proposition.}\label{2.1} 
Let $f:X \to Y$ and $g:Y \to Z$ be functions between compacta.
If $f$ and $g$ are co-existential
maps, then so is $gf$; if $gf$ is a co-existential map and $f$ is a continuous
surjection, then $g$ is a co-existential map.\\ 

\noindent
{\bf Proof.} Note that the second assertion above says more than just 
closure under terminal factors; we do not assume the co-existentiality 
of $f$.  The proof of this is immediate from the definition, so we
now consider the issue of closure under composition.   
Using \ref{2.0.5}, we have  
$p_{Y,{\cal D}}:YI\backslash {\cal D}\to Y$,     
$u:YI\backslash {\cal D}\to X$,     
$p_{Z,{\cal D}}:ZI\backslash {\cal D}\to Z$, and     
$v:ZI\backslash {\cal D}\to Y$, witnessing the co-existentiality of $f$ and $g$.
Since $g$ is co-existential, so is the ultracopower map 
$gI\backslash {\cal D}$.
Let $r:W \to ZI\backslash {\cal D}$ and        
$s:W \to YI\backslash {\cal D}$  witness the fact. 
(So $r$ is co-elementary and $s$ is a continuous 
surjection.)
Then, because compositions of co-elementary maps are co-elementary, we
have a witness to the conclusion that $gf$ is a co-existential map. $\dashv$\\

Before proceeding, some notation is in order.  If $\langle X_i:i\in I\rangle$
is a family of compacta and $\cal D$ is an ultrafilter on $I$, then, as
noted before, $\sum_{\cal D}X_i$ may be viewed as the space of maximal
filters of the normal disjunctive lattice $\prod_{\cal D}F(X_i)$.  Let
$A_i \subseteq X_i$ for
each $i \in I$.  Then we extend the ``sharp'' notation mentioned above by
defining $(\prod_{\cal D}A_i)^{\sharp}:= \{m \in \sum_{\cal D}X_i:
\prod_{\cal D}A_i \;\mbox{contains a member of}\;m\}$.  
$(\prod_{\cal D}A_i)^{\sharp}$ is closed (resp., open) in the ultracopower
if and only if $\{i\in I: A_i\;\mbox{is closed (resp., open) in}\;X_i\} 
\in {\cal D}$.  Of course, if each
$A_i$ is closed in $X_i$, and therefore a compactum in its own right,
there is a natural homeomorphism between $\sum_{\cal D}A_i$ and the
subspace $(\prod_{\cal D}A_i)^{\sharp}$.  \\     

The following lemma tells us what happens to pre-images of open sets
under codiagonal maps.  We use an overline to indicate topological
closure.\\

\subsection{Lemma.}\label{2.2}  Let $U$ be an open subset of the compactum
$X$. Then $p_{X,{\cal D}}^{-1}[U] = 
\bigcup\{(\overline{V}^I/{\cal D})^{\sharp}:
\overline{V}\subseteq U, V\;\mbox{open in}\;X\} \subseteq 
(U^I/{\cal D})^{\sharp}$.\\

\noindent
{\bf Proof.} $p_{X,{\cal D}}(m) = x$ if and only if for each open neighborhood
$W$ of $x$, $m \in (W^I/{\cal D})^{\sharp}$.  Suppose $m \in 
p_{X,{\cal D}}^{-1}[U]$, say $p_{X,{\cal D}}(m) = x$.  Let $V$ be any open set
with $x \in V \subseteq \overline{V}\subseteq U$.  Then 
$m \in (V^I/{\cal D})^{\sharp}\subseteq (\overline{V}^I/{\cal D})^{\sharp}$.

For the reverse inclusion, suppose $x = p_{X,{\cal D}}(m) \notin U$.  Pick
$V$ open with $\overline{V}\subseteq U$, and let $W$ be an open neighborhood
of $x$ such that $\overline{W} \cap \overline{V} = \emptyset$.  Then
$\overline{W}^I/{\cal D} \in m$.  But $\overline{W}^I/{\cal D} \cap
\overline{V}^I/{\cal D} = (\overline{W}\cap \overline{V})^I/{\cal D} =
\emptyset$, so $\overline{V}^I/{\cal D} \notin m$; i.e.,
$m \notin (\overline{V}^I/{\cal D})^{\sharp}$.  $\dashv$\\ 

\subsection{Theorem.}\label{2.3} Let $f:X \to Y$ be a co-existential map 
between compacta.  Then
there exists a $\cup$-homomorphism $f^*$ from the subcompacta of $Y$ to
the subcompacta of $X$ such that for each subcompactum $K$ of $Y$:
$(i)$ $f[f^*(K)] = K$; $(ii)$ $f^{-1}[U] \subseteq f^*(K)$ whenever
$U$ is open and $U \subseteq K$; and $(iii)$ the restriction $f|f^*(K)$ is a
co-existential map from $f^*(K)$ to $K$.\\

\noindent
{\bf Proof.} Let $g:YI\backslash {\cal D} \to X$ witness the fact that $f$
is a co-existential map (i.e., $fg = p_{\cal D}$), and let $K$ be a 
subcompactum of $Y$.  Then $(K^I/{\cal D})^{\sharp}$ is a subcompactum of   
$YI\backslash {\cal D}$, naturally homeomorphic to $KI\backslash {\cal D}$,
and $p_{\cal D}|(K^I/{\cal D})^{\sharp}$ may be viewed as the natural
codiagonal from $KI\backslash {\cal D}$ to $K$.  Thus $f^*(K) := 
g[(K^I/{\cal D})^{\sharp}]$ is a 
subcompactum of $X$, and $f|f^*(K)$ is a co-existential map onto its image $K$.

Suppose $U$ is open in $Y$, with $U \subseteq K$.  Then $f^{-1}[U] =
g[p_{\cal D}^{-1}[U]]$.  But by \ref{2.2}, $p_{\cal D}^{-1}[U] \subseteq
(U^I/{\cal D})^{\sharp} \subseteq (K^I/{\cal D})^{\sharp}$.  Thus
$f^{-1}[U] \subseteq f^*(K)$. 

It remains to show that the mapping $f^*$ is a $\cup$-homomorphism.  If
$K_1$ and $K_2$ are subcompacta of $Y$, then $f^*(K_1 \cup K_2) =
g[((K_1 \cup K_2)^I/{\cal D})^{\sharp}] = g[(K_1^I/{\cal D})^{\sharp} \cup
(K_2^I/{\cal D})^{\sharp}] = f^*(K_1) \cup f^*(K_2)$.   
$\dashv$\\

A number of properties, not generally preserved by continuous surjections
between compacta, are now easily seen to be preserved by co-existential maps.\\

\subsection{Proposition.}\label{2.4}  The following properties are preserved
by co-existential maps: $(i)$ being infinite; 
$(ii)$ being disconnected; $(iii)$ being a Boolean space;  
$(iv)$ being an indecomposable continuum; and $(v)$ being a hereditarily
indecomposable continuum.\\

\noindent
{\bf Proof:} Let $f,X,Y,g,{\cal D}$ and $K$ be as in the proof of \ref{2.3}.
If $Y$ is finite (resp., connected), then so is $YI\backslash {\cal D}$, and 
hence $X$.  If $Y$ is not Boolean, then we may choose $K$ to be an infinite
subcontinuum of $Y$.  This forces $f^*(K)$ to be an infinite subcontinuum 
of $X$.

Suppose $X$ is an indecomposable continuum (so there is no way to write
$X$ as the union of two proper subcontinua).  Then $Y$ is a continuum.   
Suppose $Y = K_1\cup K_2$, where each $K_n$ is a proper subcontinuum.
Then each $f^*(K_n)$ is also a proper subcontinuum of $X$, and $f^*(K_1)
\cup f^*(K_2) = f^*(K_1\cup K_2) = f^*(Y) = f^{-1}[Y] = X$.  This  
contradicts the assumption that $X$ is indecomposable.

Suppose $X$ is a hereditarily indecomposable continuum (so no subcontinuum
of $X$ is decomposable).  If $K$ is a decomposable subcontinuum of $Y$,
then $f^*(K)$ is a decomposable subcontinuum of $X$, since $f|f^*(K)$ is a 
co-existential map onto $K$.  (Use the result of the last paragraph.) 
$\dashv$\\ 

We can actually improve on \ref{2.4}$(iii)$. Noting that being a Boolean
space means being of covering dimension zero, it is tempting to conjecture
that co-existential maps cannot raise dimension.  This turns out to be the case,
but it seems we need more than just \ref{2.3} for a proof.  For us,
the handiest version of ``being of dimension $\leq n$'' 
is due to E. Hemmingsen (see \cite{E}):  A normal space
$Y$ is of covering dimension $\leq n$, $n < \omega$, ($\mbox{dim}(Y)
\leq n$) if whenever $B_1,...,B_{n+2}$ is a family of closed subsets
of $Y$ and $B_1\cap ...\cap B_{n+2} = \emptyset$, then there is a family
$F_1,...,F_{n+2}$ of closed subsets of $Y$ such that $F_m \supseteq B_m$  
for each $0 \leq m \leq n+2$, $F_1\cap ...\cap F_{n+2} = \emptyset$, and
$F_1 \cup ...\cup F_{n+2} = Y$.

\subsection{Theorem.}\label{2.4.5} Suppose $f: X \to Y$ is a 
co-existential map between compacta, $n < \omega$.  If $\mbox{dim}(X) \leq n$,
then $\mbox{dim}(Y) \leq n$.\\

\noindent
{\bf Proof.} Let $B_1,...,B_{n+2} \subseteq Y$ be given, as per the
characterization given above.       
Let $g:YI\backslash {\cal D} \to X$ witness the fact that $f$
is a co-existential map, and set $A_m := f^{-1}[B_m]$. $1 \leq m
\leq n+2$.  Each $A_m$ is closed in $X$, and $A_1\cap ...\cap A_{n+2} =
\emptyset$. 
Since $\mbox{dim}(X) \leq n$, we find appropriate closed sets $E_m 
\supseteq A_m$, $1\leq m \leq n+2$, witnessing the fact.  Pull these
sets back to the ultrapower via $g$.  We get a closed cover whose intersection
is empty.  Using compactness, we can obtain basic closed sets
$\sum_{\cal D}F_{m,i} \supseteq g^{-1}[E_m] \supseteq B_mI\backslash {\cal D}$, 
$1 \leq m \leq n+2$, such
that $\bigcap_{m=1}^{n+2}\sum_{\cal D}F_{m,i} = \emptyset$.  By elementary
ultraproduct considerations, we then have $\{i\in I: F_{m,i} \supseteq
B_m\;\mbox{for all}\;1 \leq m \leq n+2\;\mbox{and}\;
\bigcap_{m=1}^{n+2}F_{m,i} = \emptyset\;\mbox{and}\;
\bigcup_{m=1}^{n+2}F_{m,i} = Y\} \in {\cal D}$.  For any $i$ in this set,
then, we have witnesses for the fact that $\mbox{dim}(Y) \leq n$.
$\dashv$\\

A continuous map between compacta is called {\bf weakly confluent} in
the literature (see \cite{Nad}) if every subcontinuum in the range is
the image of a subcontinuum in the domain.  It follows from \ref{2.3}, plus
the fact that connectedness is preserved by co-elementary equivalence
(and continuous images, of course), that co-existential maps are weakly 
confluent.
There are many related notions discussed in \cite{Nad}, and weak confluency
is the most general.  For example, {\bf confluency} in a mapping says that each
component of the pre-image of a subcontinuum maps onto the subcontinuum.
We do not know whether there is a direct relation between co-existential maps 
and confluent maps, but suspect not.  

The strongest notion along these lines is monotonicity.
A continuous map is called {\bf monotone} if the pre-images of points
under that map are connected sets.  In \cite{Ban8} it is shown that 
co-elementary maps onto locally connected compacta are monotone; actually
this remains true for co-existential maps.\\

\subsection{Theorem.}\label{2.5}  Let $f:X \to Y$ be a co-existential map 
between compacta.  If $Y$ is locally connected, then $f$ is monotone.\\

\noindent
{\bf Proof.} Let $y\in Y$ be given, and let $\cal U$ be a neighborhood basis
for $y$ consisting of connected open sets.  For each subcontinuum 
$\overline{U}$, $U \in {\cal U}$, let $f^*(\overline{U})$  
be the subcontinuum of $X$ guaranteed to exist by \ref{2.3} 
($f^*(\overline{U})$ being connected by virtue of being a continuous image
of an ultracopower of $\overline{U}$).  Since $f^{-1}[U] \subseteq 
f^*(\overline{U})$,
we have $f^{-1}[\{y\}] \subseteq f^*(\overline{U})$.  
Now $\cal U$ is a directed set under
the ordering of reverse inclusion.  Moreover, by \ref{2.3}, $f^*$ is a
$\cup$-homomorphism, hence order preserving.  Thus the family ${\cal V}:=  
\{f^*(\overline{U}): U \in {\cal U}\}$ is also directed under
reverse inclusion.  Since $\bigcap{\cal U} = \{y\}$, it follows that
$\bigcap {\cal V} = f^{-1}[\{y\}]$.  Directed intersections of 
subcontinua are subcontinua (see \cite{Wil}); consequently $f^{-1}[\{y\}]$
is connected. $\dashv$\\  
 
The next result, an easy application of our methods, deals with isolated 
points.\\

\subsection{Proposition.}\label{2.6} Let $f:X \to Y$ be a co-existential map 
between compacta.  

$(i)$ If $y$ is an isolated point of $Y$, then $f^{-1}[\{y\}]$ is an open
singleton in $X$.

$(ii)$ If $x$ is an isolated point of $X$ and $Y$ is locally connected,
then $f(x)$ is an isolated point of $Y$.\\

\noindent
{\bf Proof.} To prove $(i)$, let $y \in Y$ be isolated.
Letting $K$ be $\{y\}$, the corresponding set $f^*(K)$ guaranteed by \ref{2.3}
must be a singleton open set, and must therefore be the pre-image under $f$
of $y$. 

To prove $(ii)$, note that $f^{-1}[\{f(x)\}]$ is connected by \ref{2.5}.  
Since $\{x\}$
is open in $X$, we have $f^{-1}[\{f(x)\}] = \{x\}$.  $\{f(x)\}$ is 
therefore open since continuous surjections between compacta are quotient
maps. $\dashv$\\ 

The question naturally arises as to the relationship between co-existential 
maps and co-elementary maps; as yet we have no examples of co-existential maps 
that are not co-elementary.  To remedy the situation, let us define an
{\bf arc} to be any homeomorphic copy of the closed unit interval in the
real line.  Arcs are of central importance in the study of compacta
(especially continua); they can be characterized as being metrizable
continua with exactly two noncut points (a classic result of R. L. Moore).  
In \cite{Ban3} it is shown (using another classic result of Moore)
that any Peano continuum (i.e., locally connected metrizable
continuum) is an arc, as long as it is co-elementarily equivalent to an arc.
It follows immediately that co-elementary images of arcs are arcs.  We can
even do a little better.\\

\subsection{Proposition.}\label{2.7} Co-existential images of arcs are arcs.\\

\noindent
{\bf Proof.} Let $f:X \to Y$ be a co-existential map between continua, where
$X$ is an arc.  $X$ is locally connected; hence so is $Y$.  By \ref{2.5},
then, $f$ is monotone.  Now $Y$ is an infinite metrizable continuum.  
So the monotonicity
of $f$, together with the fact that $X$ is an arc, tells us that $Y$ has
exactly two noncut points.  Hence $Y$ is an arc. $\dashv$\\ 

The main result of \cite{Ban8} is that any monotone continuous surjection
between two arcs is a co-elementary map.  This gives us the following.\\

\subsection{Theorem.}\label{2.8} If $f:X \to Y$ is a co-existential map 
between compacta, and if $X$ is an arc, then $f$ is co-elementary.\\

\noindent
{\bf Proof.} $Y$ is an arc by \ref{2.7}, $f$ is monotone by \ref{2.5}.
$f$ is therefore co-elementary by the main result (Proposition 2.7)
of \cite{Ban8}. $\dashv$\\
 
\subsection{Proposition.}\label{2.9} If $f:X \to Y$ is a co-elementary
map between compacta, and if $X$ is an arc, then for any subcontinuum
$C$ of $X$, $f|C$ is co-elementary onto its image if and only if $C$
and $f[C]$ have the same cardinality.\\

\noindent
{\bf Proof.} In the case $C$ is a singleton, there is no problem.  So
suppose $C$ is nondegenerate.  Then $C$ is an arc.  If $f[C]$ is a singleton,
then co-elementarity fails of course.  Otherwise, $f[C]$ is an arc too.
$f|C:C \to f[C]$ is clearly a monotone continuous surjection, and is
therefore co-elementary by Proposition 2.7 of \cite{Ban8}. $\dashv$\\

Co-elementary equivalence preserves covering dimension by results
of \cite{Ban2}.  Thus any example of a co-existential map that changes 
dimension is an example of a co-existential map that is not co-elementary.  
(Of course,
by \ref{2.4.5}, a co-existential map that changes dimension must necessarily 
lower it.)\\

\subsection{Example.}\label{2.10} We construct a dimension-lowering
co-existential map $f: X \to Y$ between locally connected metrizable continua 
as follows.

Set $X := ([0,1/2]\times\{0\})\cup [1/2,1]^2$, $Y := [0,1/2]\times \{0\}$,
$Z := [0,1]\times \{0\}$, all subsets of the euclidean plane.  $Y$ and 
$Z$ are arcs; $X$ is a ``kite-with-tail.''  $Y$ is a subcontinuum of $X$
(resp., $Z$),
and we simply let $f$ (resp., $g$) retract $X$ (resp., $Z$) onto $Y$, 
collapsing the square (resp., the interval $[1/2,1]\times\{0\}$) to the
point $\langle 1/2,0\rangle$.  Since $Y$ and $Z$ are arcs, and $g$ is
a monotone continuous surjection, we know that $g$ is co-elementary by
Proposition 2.7 of \cite{Ban8}.  We now let $h$ map $Z$ to $X$ by leaving
each element of $[0,1/2]\times \{0\}$ fixed and taking $[1/2,1]\times\{0\}$
continuously onto $[1/2,1]^2$ in such a way that $\langle 1/2,0\rangle$
remains fixed.  Clearly $fh = g$, and $f$ is therefore a co-existential map.\\

\section{An Analog of the L\"{o}wenheim-Skolem Theorem (Sharper
Version).}\label{3}
Of the many assertions that lie under the rubric ``L\"{o}wenheim-Skolem,''
the one with the most ``algebraic'' phrasing takes on the form of a
factorization theorem:  Let $\cal L$ be a first-order lexicon, with
$f:A \to B$ an embedding of $\cal L$-structures.  If $\kappa$ is an infinite
cardinal number such that $|A|+|{\cal L}| \leq \kappa \leq |B|$, then there
exists an $\cal L$-structure $C$ and embeddings $g: A \to C$, $h:C \to B$
such that $|C| = \kappa$, $h$ is an elementary embedding and $f=hg$.  

For a dualized version of this in the compact Hausdorff setting, we must
eliminate reference to a first-order lexicon, as well as decide what cardinal
invariant of compacta is to take the place of the underlying-set cardinality.
If we take a clue from Stone duality, noting that the weight of a Boolean
space equals the cardinality of its clopen algebra, then we must decide upon
weight as the
invariant of choice.  (There are, of course, other reasons to choose weight,
but the one given above is the most accessible.)  We denote the weight of
$X$ by $w(X)$.
 
\subsection{Theorem.}\label{3.1} Let $f:X \to Y$ be a continuous surjection
between compacta, with $\kappa$ an infinite cardinal such that
$w(Y) \leq \kappa \leq w(X)$.  Then there is a compactum $Z$ and continuous
surjections $g:X \to Z$, $h:Z \to Y$ such that $w(Z) = \kappa$, $g$ is a
co-elementary map, and $f = hg$.\\

\noindent
{\bf Proof.} By way of a preliminary comment, this is a sharper version
of the L\"{o}wenheim-Skolem theorem proved in \cite{Ban4}.  In that
result we used the model theory of lattices to obtain $Z$ such that
$w(Z) \leq \kappa$.  There is a bit more work involved in making sure that
$w(Z)$ may be any prescribed cardinal in the interval of possibilities. 

There is no loss of generality in assuming that $Y$ is infinite; so
first assume that $\kappa = w(Y)$, and let ${\cal B} \subseteq F(Y)$
be a lattice base of cardinality $\kappa$.  Let $\varphi :=
f^F|{\cal B}$.  Then we may treat $S(\cal B)$ as $Y$ and $\varphi^S$ as
$f$.  By the usual L\"{o}wenheim-Skolem theorem, there is a normal
disjunctive lattice $A$ and lattice embeddings $\psi :{\cal B}\to A$,
$\theta : A \to F(X)$ such that $|A| = \kappa$, $\theta$ is an elementary
embedding, and $\varphi = \theta \psi$. Set $Z := S(A)$, $g := \theta^S$,
and $h := \psi^S$.  Then $g$ is a co-elementary map, $w(Z) \leq |A| =
\kappa$, and $hg = (\theta \psi)^S = \varphi^S = f$.  (Note: This much
was proved in \cite{Ban3}.)  Now $h$ maps $Z$ onto $Y$; and, since
continuous surjections between compacta cannot raise weight, we infer that
$w(Z) = \kappa$.

The rest of the argument is not model-theoretic at all.  Suppose we could
factor $f$ into continuous surjections $g:X \to Z$, $h:Z \to Y$, where
$w(Z) = \kappa$.  Then we could apply the argument in the last paragraph
to factor $g$, and we would be done.

To effect this factorization ($f = hg$), we use the Gel'fand-Na\u{\i}mark
duality theorem between {\bf CH} and the category {\bf CBA} of 
commutative (Banach) $B^*$-algebras and nonexpansive linear maps. 
(Recall \cite{Sim} that a $B^*$-algebra
is a complex Banach algebra with a unary operation $(\;)^*$ satisfying:
$(a+b)^* = a^* + b^*$; $(ab)^* = b^*a^*$; $(\lambda a)^* =
\overline{\lambda}a^*$ ($\lambda$ is a complex scalar, $\overline{\lambda}$
is the complex conjugate of $\lambda$); $\|a^*\| = \|a\|$; and $a^{**} = a$.)
For $X \in \mbox{\bf CH}$, $C(X)$ is the $B^*$-algebra of continuous
complex-valued functions on $X$; the norm is the supremum norm, and the
involution
$(\;)^*$ is defined pointwise by $\varphi^*(x) := \overline {\varphi (x)}$.
The other half of the duality is the maximal ideal space construction $M(\;)$,
restricted to objects in {\bf CBA}.

Claim 1: $w(X) = d(C(X))$, where $d$ is the {\bf density}, the least cardinality
of a dense subset of a topological space.  To see this, let $\cal U$ be the
collection of open disks with rational centers and rational radii in the 
complex plane, and suppose $\Phi \subseteq
C(X)$ is a dense subset.  Set ${\cal V}:= \{\varphi^{-1}[U]: \varphi \in
\Phi, U \in {\cal U}\}$.  $\cal V$ is an open base for $X$.  Indeed, if
$W$ is an open neighborhood of $x$ in $X$, let $\psi :X \to [0,1]$ take
$x$ to 0 and $X \setminus W$ to 1.  Let $\varphi \in \Phi$ be such that
$\|\varphi - \psi\| < 1/4$.  Then $\|\varphi (0)\| < 1/4$ and
$\|\varphi (y)\| > 3/4$ for $y \in X\setminus W$; hence if $U \in {\cal U}$ is
the open disk of radius 1/2 centered at the origin, then $\varphi^{-1}[U]$
is a set in $\cal V$ containing $x$ and contained in $W$.  Since
$|{\cal V}| \leq |\Phi|$, we have $w(X) \leq d(C(X))$.  Now let $\cal V$
be an open base for $X$.  By Weierstrass approximation, we can get a
dense $\Phi \subseteq C(X)$ of cardinality $\leq |{\cal V}|$.  Thus we
have $d(C(X)) \leq w(X)$.
 
Claim 2: There is a subset $\Phi$ of $C(X)$, consisting of maps into
the unit interval, such that: $(i)$ $|\Phi| = w(X)$; and $(ii)$ 
$\|\varphi - \psi \|
\geq 1/2$ for all $\varphi , \psi \in \Phi$.  To see this, first note
that any set $\Phi$ satisfying $(ii)$ must have cardinality 
at most $d(C(X))$ (which is $w(X)$, by Claim 1).  Now let $\Phi$ be maximal
with regard to satisfying $(ii)$, let $U$ be the open disk of radius
1/2 and centered at the origin in the complex plane, and set
${\cal V} := \{\varphi^{-1}[U]: \varphi \in \Phi \}$.  We need to show
$\cal V$ is an open base for $X$.  Assuming the contrary, there exist
$x \in W \subseteq X$, $W$ open in $X$, such that for all $\varphi \in
\Phi$, if $\varphi (x) \in U$ (so $\varphi (x) \in [0, 1/2)$), then there is
some $y \notin W$ with $\varphi (y) \in [0,1/2)$.  Let $\psi: X \to [0,1]$ take
$x$ to 0 and $X \setminus W$ to 1.  If $\varphi \in \Phi$ and 
$\varphi (x) \in [1/2,1]$, then $\|\psi - \varphi\| \geq 
|\psi (x) - \varphi (x)|
\geq 1/2$.  If $\varphi (x) \in [0,1/2)$, let $y \in X \setminus W$ be such that
$\varphi (y) \in [0,1/2)$.  Then $\|\psi - \varphi\| \geq 
|\psi (y) - \varphi (y)|
\geq 1/2$.  Thus $\Phi \cup \{\psi\}$ properly contains $\Phi$ and satisfies
$(ii)$, contradicting the maximality of $\Phi$.  Since $\cal V$ is an
open base for $X$, and $|{\cal V}| \leq \Phi$, we have $w(X) \leq |\Phi|$.
Thus $|\Phi| = w(X)$.

Now we have an embedding $f^C: C(Y) \to C(X)$.  Let $\kappa$ be 
such that $w(Y) \leq \kappa \leq w(X)$, and let $\Phi \subseteq C(X)$
have cardinality $\kappa$ and satisfy $(ii)$ in Claim 2.  Let $\Psi
\subseteq C(Y)$ be dense of cardinality $w(Y)$.  Then $\Phi \cup
f^C[\Psi]$ is a subset of $C(X)$ of cardinality $\kappa$; closure of this set
under the Banach algebra operations with scalar multiplication restricted to
the rational complex numbers then results in a subring $A$, also of cardinality 
$\kappa$.
The topological closure $\overline{A}$ of $A$ is therefore a $B^*$-subalgebra
of $C(X)$.  Because $|A| = \kappa$, $d(\overline{A}) \leq \kappa$.
Because $\Phi \subseteq \overline{A}$, $d(\overline{A}) \geq \kappa$.
Because $\Psi$ is dense in $C(Y)$, $f^C[C(Y)] \subseteq \overline{A}$.
Thus, by applying the maximal ideal functor $M(\;)$ 
(i.e., $Z:= M(\overline{A})$), we get the factorization we want.  This 
completes the proof.  $\dashv$\\

Let us return, for the moment, to the model-theoretic setting; for 
simplicity, assume the underlying lexicon $\cal L$ is countable.
A nice application of the L\"{o}wenheim-Skolem factorization theorem
mentioned in the lead paragraph of this section is the following.\\

\subsection{Proposition.}\label{3.2} 
Let $f:A \to B$ be an embedding
between infinite $\cal L$-structures, with $\kappa \leq |A|$ an infinite
cardinal.  Suppose that for each $\cal L$-structure $C$ of cardinality
$\kappa$, and each 
elementary embedding $g:C \to A$, the composition $fg$ is elementary
(resp., existential).  Then $f$ is elementary (resp., existential).\\

\noindent
{\bf Proof.} Let $\varphi (x_1,...,x_n)$ be any (existential)
formula, with $\langle a_1,...,a_n\rangle$ an $n$-tuple from $A$.  By
L\"{o}wenheim-Skolem, there is an $\cal L$-structure $C$ of cardinality
$\kappa$ and an elementary embedding $g:C \to A$ such that $g[C]$ contains
each $a_m$.  The desired conclusion is now immediate. $\dashv$\\

The proof above, while extremely simple, suffers the
major failing of being nonportable.  The Proposition has an obvious
restatement in the language of compacta, and this proof sheds very little
light on how to establish the dualized version.  There is, however, an  
argument,
based solely on ultraproduct considerations (inspired by C. C. Chang's
ultraproduct proof of the compactness theorem, see \cite{CK,Ekl}), which
can be easily transported to the topological context.  Its drawback is 
that it is relatively cumbersome.  We dualize that argument in the following.

\subsection{Theorem.}\label{3.3}
Let $f:X \to Y$ be a continuous surjection between infinite compacta,
with $\kappa \leq w(Y)$ an infinite
cardinal.  Suppose that for each compactum $Z$ of weight 
$\kappa$, and each 
co-elementary map $g:Y \to Z$, the composition $gf$ is co-elementary
(resp., co-existential).  Then $f$ is co-elementary (resp., co-existential).\\

\noindent
{\bf Proof.} In model theory, the obstruction to being able to show,
with mapping diagrams, that the class of elementary (resp., existential)
embeddings is closed under terminal factors, resides in a failure of
surjectivity: one cannot carry out a successful diagram chase. On the
topological side, the obstruction to closure under initial factors
resides in a failure of injectivity.  Both these failures can be 
remedied somewhat with the use of ultra(co)products.

We establish our result for the co-elementary case; the co-existential case
is handled in like fashion.   
First let $\Delta$ be the set of finite subsets of $Y$.  For
each $\delta \in \Delta$ there is a continuous mapping $r_{\delta}$ from 
$Y$ to the closed unit
interval, such that the restriction $r_{\delta}|\delta$ is one-one.  
Let $W_{\delta} :=
r_{\delta}[Y]$.  By \ref{3.1}, there is a compactum $Z_{\delta}$ of weight 
$\kappa$,
and continuous surjections $g_{\delta}:Y \to Z_{\delta}$, 
$t_{\delta}:Z_{\delta} \to W_{\delta}$, such
that $g_{\delta}$ is co-elementary and $r_{\delta} = t_{\delta}g_{\delta}$.  
So each $g_{\delta}$ is a
co-elementary map onto a compactum of weight $\kappa$, and 
$g_{\delta}|\delta$ is one-one.

Now, with the aid of \ref{2.0.5}, we build a diagram $\mbox{\bf D}_{\delta}$
that witnesses the co-elementarity above,
for each $\delta \in \Delta$.  More precisely, we have a single ultrafilter
$\cal D$ on a set $I$, and homeomorphisms 
$h_{\delta}:XI\backslash {\cal D} \to
Z_{\delta}I\backslash {\cal D}$, $k_{\delta}:YI\backslash {\cal D} \to 
Z_{\delta}I\backslash {\cal D}$ witnessing the co-elementarity of 
$g_{\delta}f$ and
$g_{\delta}$ respectively.  Letting $p_{\delta}:Z_{\delta}I\backslash {\cal D} 
\to Z_{\delta}$,
$p_X:XI\backslash {\cal D} \to X$, and $p_Y:YI\backslash {\cal D} \to Y$
be the canonical codiagonal maps, we then have $p_{\delta}k_{\delta} = 
g_{\delta}p_Y$ and
$p_{\delta}h_{\delta} = g_{\delta}fp_X$.

In general we cannot expect $fp_X = p_Yk_{\delta}^{-1}h_{\delta}$; 
however it is true
that $g_{\delta}fp_X = g_{\delta}p_Yk_{\delta}^{-1}h_{\delta}$.  We now
proceed to take an ``ultracoproduct'' of the diagrams $\mbox{\bf D}_{\delta}$.

For each $\delta \in \Delta$, let $\hat{\delta} := \{\gamma \in \Delta:
\delta \subseteq \gamma \}$.  Then the set $\{\hat{\delta}: \delta \in
\Delta\}$ clearly satisfies the finite intersection property, and hence
extends to an ultrafilter $\cal H$ on $\Delta$.  
Form the ``$\cal H$-ultracoproduct'' diagram {\bf D} in the obvious way.
Then we have the canonical codiagonal maps $u_X:X\Delta\backslash {\cal H}
\to X$, $v_X:(XI\backslash {\cal D})\Delta \backslash {\cal H} \to
XI\backslash {\cal D}$, 
$u_Y:Y\Delta \backslash {\cal H}
\to Y$, and $v_Y:(YI\backslash {\cal D})\Delta \backslash {\cal H} \to
YI\backslash {\cal D}$.  Moreover, $p_Xv_X = u_X(p_X\Delta \backslash {\cal H})$
and $p_Yv_Y = u_Y(p_Y\Delta \backslash {\cal H})$ are (essentially) codiagonal
maps from iterated ultracopowers, and these ultracopowers are isomorphic
via $(\sum_{\cal H}k_{\delta})^{-1}(\sum_{\cal H}h_{\delta}) =
\sum_{\cal H}k_{\delta}^{-1}h_{\delta}$.

We need to show that $fp_Xv_X = p_Yv_Y
(\sum_{\cal H}k_{\delta})^{-1}(\sum_{\cal H}h_{\delta})$.  Suppose otherwise.  
Then we have some $x \in (XI\backslash {\cal D})\Delta \backslash {\cal H}$
with $[fp_Xv_X](x) = y_1 \neq y_2 = [p_Yv_Y
(\sum_{\cal H}k_{\delta})^{-1}(\sum_{\cal H}h_{\delta})](x)$.  
Because $\cal H$ is an ultrafilter, we have
$\{\gamma \in \Delta: 
v_Y(\sum_{\cal H}k_{\delta})^{-1}(\sum_{\cal H}h_{\delta}) = 
k_{\gamma}^{-1}h_{\gamma}v_X\} \in {\cal H}$; hence
$\{\gamma \in \Delta: y_2 =
[p_Yk_{\gamma}^{-1}h_{\gamma}v_X](x)\} \in {\cal H}$.  Now
$[g_{\gamma}p_Yk_{\gamma}^{-1}h_{\gamma}v_X](x) = [g_{\gamma}fp_Xv_X](x)$ for
all $\gamma \in \Delta$, so $\{\gamma \in \Delta: g_{\gamma}(y_2) =
g_{\gamma}(y_1)\} \in {\cal H}$.  However, 
$\{\gamma \in \Delta: g_{\gamma}(y_2) \neq 
g_{\gamma}(y_1)\} \supseteq \widehat{\{y_1,y_2\}} \in {\cal H}$.  
This contradiction completes the proof. $\dashv$\\

\section{An Analog of the Elementary Chains Theorem.}\label{4}
The simplest version of the elementary chains theorem, due jointly to A. Tarski
and R. L. Vaught, says that the union of an $\omega$-indexed elementary chain
of relational structures is an elementary extension of each of its summands.
The obvious translation into the realm of compacta,
what we refer to here as the {\bf co-elementary chains hypothesis} (CECH)
is the assertion that the inverse limit
of an $\omega$-indexed co-elementary chain of compacta is a co-elementary
cover of each of its factors. More precisely, if
$\langle X_n \stackrel{f_n}{\leftarrow} X_{n+1}: n < \omega \rangle$  
is an $\omega$-indexed family of co-elementary maps, and if $X :=
\displaystyle \lim_{\leftarrow}X_n$ is the inverse limit space with natural
connecting maps $g_n:X \to X_n$, $n < \omega$, then each $g_n$ is a 
co-elementary map. (Recall that $X$ is defined to be the subspace
$\{ \langle x_0,x_1,\dots \rangle \in  \prod_{n<\omega}X_n:
x_n = f_{n+1}(x_{n+1})\; \mbox{for all} \; n < \omega\}$, and $g_n$ is 
projection to the $n\:$th factor, restricted to this subspace.)\\  

The way one normally goes about proving the usual elementary chains theorem
is to use induction on the complexity of formulas.  We know of no other
way; in particular, we know of no way to prove the result using the
ultrapower theorem and chains of isomorphisms of ultrapowers.  This  
perhaps speaks to an inherent lack of ``fine structure'' in 
ultraproduct methods and definitely
presents difficulties when one tries to find a direct proof of the 
CECH.  While the CECH is true, the only proof we know of uses 
Gel'fand-Na\u{\i}mark duality and Banach model theory (see below).

The following weak version of the CECH {\it can} be proved by direct
methods (avoiding Banach model theory, at any rate), and is worth while
exploring.  The difference between the weak and strong versions of
the CECH may be compared to the difference between the amalgamation property 
for elementary
embeddings (resp., co-elementary maps) and the existence of elementary
pushouts (resp., co-elementary pullbacks).\\   

\subsection{Proposition.}\label{4.1} 
Let $\langle X_n \stackrel{f_n}{\leftarrow} X_{n+1}: n < \omega \rangle$  
be a co-elementary chain of compacta.  Then there exists a compactum $X$
and co-elementary maps $g_n:X \to X_n$, $n<\omega$, such that for all
$n<\omega$, $g_n = f_ng_{n+1}$.\\

\noindent
{\bf Proof.} 
By \ref{2.0.5}, there is an ultrafilter $\cal D$ on a set $I$ and
homeomorphisms
$X_nI\backslash {\cal D} \stackrel{h_n}{\leftarrow}
X_{n+1}I\backslash {\cal D}$, $n < \omega$, such that for each $n$,
$p_{X_n,{\cal D}}h_n = f_np_{X_{n+1},{\cal D}}$.  Let $X$
be the inverse limit of this chain, with natural 
maps $k_n:X \to
X_nI\backslash {\cal D}$.  Since each $h_n$ is a homeomorphism, so
is each $k_n$, and we set $g_n := p_{X_n,{\cal D}}k_n$. $\dashv$\\

Of course the weakness of our weak version of the CECH in \ref{4.1} lies
in the fact that the compactum $X$ is not generally $\displaystyle
\lim_{\leftarrow}X_n$.  This would be remedied if we could show that the
natural map $h:X \to \displaystyle \lim_{\leftarrow}X_n$, a continuous
surjection, is actually co-elementary; but that seems to be beyond our
methods.

The only way we know of to prove the CECH is to use Banach model theory.
For reasons of space, our argument is far from self-contained; the
interested reader is referred to the appropriate literature 
(e.g., \cite{HH, HHM, Hen1, Hen2, HI, Sim}) for details.  The skeleton
of the approach is this: 
The basic lexicon of Banach model theory uses the usual symbols from the 
theory of 
abelian groups, plus countably many unary operation symbols to allow for scalar 
multiplication.  (The scalar field can be either the field of real rationals
or complex rationals.)  Also there is a unary predicate symbol whose intended
interpretation is the closed unit ball in a Banach space.  (The theory is
flexible enough to allow for additional operation symbols; e.g., lattice
symbols, multiplication and involution.)  

Banach model theory does not use all first-order formulas, just the
positive-bounded ones (i.e., conjunction, disjunction, and quantification
restricted to the unit ball).  Also the relation of satisfaction is
weakened to what is called {\it approximate satisfaction}.  Once this is
all laid out, the notions of elementary equivalence and elementary
embedding make sense in the Banach setting.  Finally there is an 
analogous ultrapower theorem, where ``ultrapower'' means Banach ultrapower
and ``isomorphic'' means isometrically isomorphic.
This tells us that the Banach notion of elementary embedding is the 
exact analog of our notion of co-elementary map.  

We now come to the main bridge connecting the Banach world and the
world of compacta, namely the Gel'fand-Na\u{\i}mark duality, a genuine 
category-theoretic duality between
{\bf CH} and the category {\bf CBA} of commutative $B^*$-algebras and
nonexpansive linear maps (see the proof of \ref{3.1} above).  
Using a result called the perturbation
lemma \cite{Hen2}, the obvious analog of the classic elementary chains
theorem is proved in the Banach setting \cite{HI}.  
(Connecting maps are elementary
embeddings in the Banach sense, and the full direct limit is the direct
limit in the category {\bf CBA} when the summands happen to be commutative
$B^*$-algebras.) So the CECH is nothing but the Gel'fand-Na\u{\i}mark
dual of the Banach elementary chains theorem restricted to objects in
{\bf CBA}.  We thus have what is essentially a proof of the CECH.\\

\subsection{Theorem.}\label{4.2} In \ref{4.1}, $X$ may be taken to
be the inverse limit $\displaystyle \lim_{\leftarrow}X_n$.\\ 
 
\subsection{Remark.}\label{4.3}  The technique for proving the
elementary chains theorem allows for generalization to directed systems
of elementary embeddings.  This then gives us a corresponding generalized
CECH.\\

\section{An Analog of Robinson's Test.}\label{5}
The notion of model completeness in first-order model theory was invented
by A. Robinson, who was inspired by classical (Nineteenth-Century)
algebra; in particular, Hilbert's {\it Nullstellensatz} concerning the location
of solutions of systems of equations and inequations.  (See, e.g., \cite{Hod,
Mac}.)
One of many equivalent formulations of this notion is to say that a first-order
theory is {\bf model complete} if any embedding between two models of the
theory is elementary.
Since elementary embeddings are a great deal rarer than embeddings in
general, it is very interesting when a first-order theory is discovered to be
model complete.  (Examples include (see \cite{CK}):  dense linear orderings 
without 
endpoints; atomless Boolean algebras; algebraically closed fields; real closed
fields; real closed ordered fields.)  It is not surprising, then, that much
effort has been expended in the study of this phenomenon; especially in the
search for readily applicable tests to detect its presence.
 
Robinson's test says that a theory is model complete just in case every
embedding between models of the theory is an existential embedding.  Since 
testing
for existentiality is ostensibly easier than testing for elementarity, this
result is very important in the general study (as well as being a
strikingly elegant application of the elementary chains theorem).  

To give Robinson's test a proper phrasing in the compact Hausdorff setting,
we first define a class of compacta to be {\bf co-elementary} 
if it is closed under co-elementary equivalence and the taking
of ultracopowers.  (This corresponds precisely to a class of relational
structures being the class of models of a first-order theory; examples
include (see \cite{Ban2}): Boolean spaces without isolated points; continua;
(in)decomposable continua; compacta of covering dimension $n$, $n < \omega$;
infinite-dimensional compacta.   
The class of locally connected compacta is {\it not} co-elementary
\cite{Gur,Ban5}.)  A co-elementary
class is {\bf model complete} if any continuous surjection between two of its
members is a co-elementary map.  (Aside:  recalling what completeness
means in model theory, the right criterion for a co-elementary class to
be {\bf complete} is that it consist of just one co-elementary equivalence 
class.  It is not entirely trivial to show that co-elementary equivalence
classes are indeed co-elementary classes, but true nonetheless. (See 
\cite{Ban7}.  Lemma \ref{2.0.5} greatly facilitates the proof that
co-elementary equivalence classes are closed under ultracoproducts.)  
The notions of completeness and model completeness, while cognate,
are not directly related logically.  An obvious analog of the prime model
test, also obviously true in the topological setting, is that a model
complete co-elementary class {\bf K} is complete if there is some $X \in
\mbox{\bf K}$ such that every member of {\bf K} continuously surjects onto
$X$.)  It is our belief that the reasons
for studying model completeness in model theory remain just as compelling
when we consider the topological analog; the present work is just a 
beginning of a process of discovering interesting co-elementary classes
that are model complete.    

\noindent
\subsection{Theorem.}\label{5.1} A co-elementary class of compacta is
model complete if and only if every continuous surjection in the class
is a co-existential map.\\

\noindent
{\bf Proof.}  Our proof is an exact dualization of the well-known proof
of Robinson's test using elementary chains (see, e.g., \cite{Mac}); we
present it here for the sake of thematic completeness.   

Only one direction is nontrivial.  Assume that {\bf K}
is a co-elementary class with the property that every continuous surjection
in {\bf K} is a co-existential map, and let $f:X \to Y$ be a continuous 
surjection
in {\bf K}.  $f$ is co-existential, so we get a compactum $Z_1$ and continuous
surjections $g_0:Z_1 \to Y$, $h_0:Z_1 \to X$ with $g_0$ co-elementary
and $g_0 = fh_0$.  $Z_1$ is thus in {\bf K}, and we conclude that $h_0$
is co-existential.  Repeat the procedure.  We get a compactum $Z_2$ and 
continuous
surjections $g_1:Z_2 \to X$, $h_1: Z_2 \to Z_1$ with $g_1$ co-elementary
and $g_1 = h_0h_1$.  $Z_2 \in \mbox{\bf K}$, and $h_1$ is thus a co-existential
map.   
Proceeding in this way, we obtain compacta $Z_n \in \mbox{\bf K}$ and
continuous surjections $g_n,h_n$, $n < \omega$.  For $n \geq 1$,
$h_n: Z_{n+1} \to Z_n$, $g_{n+1}:Z_{n+2} \to Z_n$, $g_{n+1} =
h_nh_{n+1}$.  Each $h_n$ is co-existential, each $g_n$ is co-elementary.
The inverse limit $Z$ of the maps $h_n$, $n < \omega$, is also the inverse limit
of the co-elementary maps $g_n$, $n$ odd (resp., $n$ even).
By the CECH (\ref{4.2}), there are co-elementary maps
$u: Z \to X$ and $v: Z \to Y$ such that $v = fu$.  Since the class of
co-elementary maps is closed under terminal factors, it follows that
$f$ is co-elementary. $\dashv$\\ 

There is a sharper version of Robinson's test
(see, e.g., \cite{CK}).
One assumes that the theory has no finite models; and that, for some infinite
cardinal $\kappa \geq |{\cal L}|$, any embedding between models of
cardinality $\kappa$ is an existential embedding.  Then the theory is model 
complete.  The following is the topological analog of this fact.

\subsection{Theorem.}\label{5.2}
Let {\bf K} be a co-elementary class containing no finite compacta.
Then {\bf K} is model complete if (and only if) there is some infinite 
cardinal $\kappa$ such
that each continuous surjection between compacta of weight $\kappa$ in 
{\bf K} is a co-existential map.\\ 

\noindent
{\bf Proof.} Let $f:X \to Y$ be a continuous surjection in {\bf K}.
By \ref{5.1}, it suffices to show that $f$ is co-existential.  
Assume first that $\kappa \leq w(Y)$.  By \ref{3.3}, it suffices to
show that $gf$ is co-existential whenever $g:Y \to Z$ is co-elementary
and $w(Z) = \kappa$.  So let $g$ be given.    
By \ref{3.1}, we can obtain a factorization $u:X \to W$, $v:W \to Z$,
where $u$ is co-elementary, $w(W) = w(Z) = \kappa$, and $gf = vu$.
Of course both $W$ and $Z$ are in {\bf K}, so $v$ is co-existential.
Then $vu$ is co-existential by \ref{2.1}.

If $\kappa > w(Y)$, then find an ultrafilter $\cal D$ on a set $I$ such that
$w(YI\backslash {\cal D}) \geq \kappa$ (see \cite{Ban2}).  
By the argument above, we conclude
that $fI\backslash {\cal D}$ is co-existential.  Immediately we infer
(see \ref{2.1}) that $f$ is co-existential. $\dashv$\\

\section{Co-inductive Co-elementary Classes.}\label{6}
The usual elementary chains theorem shows that any model complete theory
is {\bf inductive}; i.e., closed under chain unions (indeed, direct limits
of directed systems of embeddings).  Being inductive, by the 
Chang-\L o\'{s}-Suszko theorem, is equivalent to being AE axiomatizable. 
A co-elementary class {\bf K} is {\bf co-inductive} if it is closed under
inverse limits of directed systems of continuous surjections.  
Because of \ref{4.2},
we know that model complete co-elementary classes are co-inductive.

The following notion dualizes that of being existentially closed relative to
an elementary class of relational structures.
Let {\bf K} be a co-elementary class, with $X \in \mbox{\bf K}$.  $X$
is {\bf co-existentially closed} in {\bf K} if whenever $f: Y \to X$ is a 
continuous surjection and $Y \in \mbox{\bf K}$, then $f$ is a co-existential
map.

\subsection{Theorem.}\label{6.1} Let {\bf K} be a co-inductive
co-elementary class, with $X \in \mbox{\bf K}$ infinite.  Then $X$ is 
the continuous image of some $Y$ that is co-existentially closed in {\bf K}, 
such that $w(Y) = w(X)$.\\

\noindent
{\bf Proof.} This is the dualization of a well-known property of 
inductive elementary classes, so we supply just the transition steps.

The class of normal disjunctive lattices whose maximal spectra lie in
{\bf K} is suitably denoted $S^{-1}[\mbox{\bf K}]$.  The class {\bf NDL}
of normal disjunctive lattices is itself an inductive elementary class.
Since $S(\;)$ converts direct limits to inverse limits, we conclude that
$S^{-1}[\mbox{\bf K}]$ is an inductive elementary class whenever
{\bf K} is a co-inductive co-elementary class.  
Suppose $X \in \mbox{\bf K}$ is infinite, and pick a lattice
base $\cal A$ for $X$ of infinite cardinality $w(X)$.  
Then ${\cal A} \in
S^{-1}[\mbox{\bf K}]$, and there is an embedding
$f :{\cal A} \to B$ for some existentially closed $B \in S^{-1}[\mbox{\bf K}]$
such that $|B| = |{\cal A}|$.  (The proof of this is elementary model
theory, involving unions of chains of embeddings (see \cite{CK,Hod}).)  
Set $Y := S(B)$.
Then $f^S$ is (essentially) a continuous surjection from $Y$ to $X$, so
$w(X) \leq w(Y)$.  But $w(Y) \leq |B| = |{\cal A}| = w(X)$, so the weights
are equal.  Finally, if $g: Z \to Y$ is a continuous surjection, $Z \in
\mbox{\bf K}$, then $g^F: F(Y) \to F(Z)$ is an embedding in   
$S^{-1}[\mbox{\bf K}]$.  Let $u: B \to F(Y)$ be the obvious separative
embedding.  Since $B$ is existentially closed, $g^Fu$ is an existential
embedding.
Since $S(\;)$ converts existential embeddings to co-existential maps (and 
separative embeddings
to homeomorphisms), we infer that $g$ is a co-existential map. $\dashv$\\  

By way of a linguistic aside, suppose ``blob'' is a noun given to name a 
class of compacta that is a co-inductive co-elementary class.  Then we
may abbreviate ``$X$ is co-existentially closed in the class of blobs''
as ``$X$ is a co-existentially closed blob.''  Note, however, that this 
does not mean that $X$ is a blob that is co-existentially closed in the
class of compacta, any more than ``$G$ is a free abelian group'' means
that $G$ is a free group that happens to be abelian.  (Or, more prosaically,
any more than ``George is a small elephant'' means that George is a small
entity that happens to be an elephant.)  (W. V. O. Quine \cite{Q} applies 
the word {\it syncategorematic} to such adjectives as ``co-existentially
closed,'' ``free,'' and ``small'' used in this way.)  Because of 
\ref{6.1}, we know that there exist co-existentially closed compacta in
all infinite weights.  The same can be said for the existence of
co-existentially closed continua, since the class of continua is
co-inductive co-elementary.  Characterizing the co-existentially closed
blobs, then, is a new problem each time we change the meaning of the word
``blob.''  When this meaning is least restrictive, the problem turns out 
to have an easy solution.\\

\subsection{Proposition.}\label{6.1.1}  The co-existentially closed
compacta are precisely the Boolean spaces without isolated points.\\

\noindent
{\bf Proof.}  We call a space without isolated points {\bf self-dense}.
Suppose $X$ is a co-existentially closed compactum.  Then $X$ is a
co-existential image of a Stone-\v{C}ech compactification $\beta (I)$
for a suitably large discrete space $I$.  By \ref{2.4}$(iii)$, $X$ must
therefore be Boolean.  By Stone duality, its clopen algebra $B(X)$ must
be an existentially closed Boolean algebra.  Now the class of such algebras
comprises the atomless Boolean algebras (see \cite{Mac}); hence $X$ is
a self-dense Boolean space.
 
In the other direction, suppose $X$ is a self-dense Boolean space, with
$f:Y \to X$ a continuous surjection between compacta.  Let $I$ be a discrete
space of cardinality $|Y|$, and set $Z := \beta (I) \times W$, where
$W$ is any self-dense Boolean space.  Then $Z$ is a self-dense Boolean space,
and there is a continuous surjection $g :Z \to Y$.  Let $h := fg$.  Then
$h$ is a continuous surjection between self-dense Boolean spaces; hence
$h^B :B(X) \to B(Z)$ is an embedding between atomless Boolean algebras.
The class of such algebras is model complete; therefore $h^B$ is elementary.
This says that $h$ is a co-elementary map, proving that $f$ is co-existential.
$\dashv$\\
 
The problem of identifying the co-existentially closed continua is still
open, and seems to be quite difficult.  What little we know so far, besides
the fact that they exist in all infinite weights, is the following.

\subsection{Proposition.}\label{6.1.2}  Every co-existentially closed
continuum is indecomposable.  (Consequently, every nondegenerate continuum
is a continuous image of an indecomposable continuum of the same weight.)\\

\noindent
{\bf Proof.} 
Suppose $X$ is a decomposable continuum.  Then (see \cite{Kur,Wil}) $X$
has a proper subcontinuum $K$ and a nonempty open set $U \subseteq K$.
Pick $x_0 \in X\setminus K$, and
let $Y := (X \times \{0,1\})/\!\!\sim$ be the quotient space where the only
nontrivial identification is $\langle x_0,0\rangle \sim \langle x_0,1\rangle$.
Let $f:Y \to X$ be induced by projection onto the first factor.  Then
$f$ is a continuous surjection between continua, and $f^{-1}[U]$ 
intersects both components of $f^{-1}[K]$.  By \ref{2.3}, then, $f$
cannot be a co-existential map and $X$ is therefore not a co-existentially
closed continuum.  The parenthetical assertion is now an instant corollary
of \ref{6.1} and the above.  $\dashv$\\  

In model theory, there is an important connection between model completeness 
and categoricity.  In particular, there is Lindstr\"{o}m's test (see
\cite{Hod, Mac}), which 
says that if a consistent AE theory $T$ has no finite models and is
$\kappa$-categorical for some infinite cardinal $\kappa \geq |{\cal L}(T)|$,
then $T$ is model complete.  The techniques developed here allow us to
establish a topological analog of this.  First define a 
co-elementary class {\bf K} to be $\kappa$-{\bf categorical}, where 
$\kappa$ is an infinite cardinal, if: $(a)$ there are compacta of weight
$\kappa$ in {\bf K}; and $(b)$ any two compacta of weight $\kappa$ in
{\bf K} are homeomorphic.  (Clearly, by the L\"{o}wenheim-Skolem theorem
(\ref{3.1}), if {\bf K} is a co-elementary class containing no finite
compacta, and if {\bf K} is $\kappa$-categorical for some cardinal
$\kappa$, then {\bf K} is complete.) 

\subsection{Theorem.}\label{6.1.5} 
Let {\bf K} be a co-elementary
class containing no finite compacta.  If {\bf K} is co-inductive and
$\kappa$-categorical for some infinite cardinal $\kappa$, then {\bf K} is
model complete.\\

\noindent
{\bf Proof.} By \ref{5.2}, it suffices to show that every continuous
surjection between compacta of weight $\kappa$ in {\bf K} is co-existential.
But \ref{6.1} and $\kappa$-categoricity tell us that all compacta in {\bf K}
of weight $\kappa$ are co-existentially closed in {\bf K}. $\dashv$\\

\subsection{Remarks and Questions.}\label{6.2} 
$(i)$ Concrete examples of co-elementary and co-existential maps are
hard to find.  One simple question we would like to settle is the following:
Are monotone continuous surjections between nondegenerate
locally connected continua necessarily co-existential maps?  In particular, is
the projection map from the closed unit square onto its first co\"{o}rdinate
a co-existential map?

$(ii)$ If {\bf K} is a 
co-elementary class containing
at least one infinite compactum, then
$S^{-1}[\mbox{\bf K}]$ is {\it never} model complete.  By Proposition 2.8
in \cite{Ban7}, $f^F:F(Y) \to F(X)$ is an elementary embedding just in
case $f: X \to Y$ is a homeomorphism.  The existence of codiagonal maps in
{\bf K} that are not homeomorphisms then provides embeddings in
$S^{-1}[\mbox{\bf K}]$ that are not elementary.
 
$(iii)$ Define a compactum $X$ to be {\bf $\kappa$-categorical} if:
$(a)$ the co-elementary equivalence class of $X$ 
is $\kappa$-categorical; and $(b)$ $w(X) = \kappa$.  $X$ is {\bf categorical}
if $X$ is $w(X)$-categorical.  The only categorical compacta we know of are
Boolean (the Cantor discontinuum, for example, being $\aleph_0$-categorical).
By the main result of \cite{Ban3}, arcs are ``categorical''
in a restricted sense:  one may look only among the locally connected
compacta.  Indeed, every infinite compactum is co-elementarily equivalent
to a compactum (of any prescribed infinite weight) that is not locally 
connected \cite{Ban5}, so any examples of infinite categorical compacta
must fail to be locally connected.  In fact, if $X$ is $\kappa$-categorical  
and $Y \equiv X$ has weight $\geq \kappa$, then $Y$ co-elementarily
surjects onto some compactum $Z$ of weight $\kappa$.  $Z$ is then homeomorphic
to $X$, which is not locally connected.  It follows that $Y$ is not
locally connected either.  In particular, no compactum co-elementarily 
equivalent to an $\aleph_0$-categorical compactum  can be locally connected.
For any cardinal $\kappa$ (finite too), define
a compactum $X$ to be $\kappa$-{\bf wide} if for each cardinal $\lambda 
< \kappa$, $X$ contains a family of
$\lambda$ pairwise disjoint proper subcontinua with nonempty interiors.  
Clearly an infinite locally connected compactum is $\aleph_1$-wide; and 
decomposability for continua is equivalent to being 2-wide.  
By techniques similar to
those used to show that the classes of (in)decomposable continua are
co-elementary (see \cite{Ban5,Gur}), one can prove that the class of
$\aleph_0$-wide compacta is co-elementary (but not its complement).
(Indeed, for $n<\omega$, $\sum_{\cal D}X_i$ is $n$-wide if and only if 
$\{i\in I:
X_i\;\mbox{is}\; n\mbox{-wide}\} \in {\cal D}$.)  So if
$X$ is co-elementarily equivalent to an infinite locally connected compactum,
then $X$ is $\aleph_0$-wide.  We do not think that the converse is true.  
(In a private conversation,
C. W. Henson has told us that an $\aleph_0$-categorical compactum must
fail to be $\aleph_0$-wide.  His method uses a Banach version of the
Ryll-Nardzewski theorem.)     

$(iv)$ We still lack examples of model complete co-elementary classes
that are not subclasses of the Boolean spaces.

$(v)$ We lack examples of pairs of co-elementarily equivalent Peano
continua that are not homeomorphic.

$(vi)$ What, besides \ref{6.1.2}, can we infer about co-existentially  
closed continua?  (In view of \ref{2.4.5}, for example,
we suspect that co-existentially closed continua are {\bf curves}; i.e., 
one-dimensional.)  Dare we hope for a classically topological 
classification?  Are any of the familiar examples of indecomposable
continua (e.g., pseudo-arcs, solenoids) co-existentially closed?
 
$(vii)$ Let {\bf CCH} be the class of continua (considered as a full
subcategory of {\bf CH}).  There is a single AE sentence $\gamma$ in 
the first-order language of bounded lattices such that a normal disjunctive
lattice $A$ satisfies $\gamma$ if and only if $S(A) \in \mbox{\bf CCH}$
(see \cite{Ban7}).  Let {\bf CNDL}$:= \{A \in \mbox{\bf NDL}: A \models
\gamma\} = S^{-1}[\mbox{\bf CCH}]$.  Then, as in the proof of \ref{6.1},
each existentially closed member of {\bf CNDL} gives rise to a co-existentially
closed member of {\bf CCH}.  Is this assignment surjective?  Is the class
of existentially closed members of {\bf CNDL} an elementary class?
(This deals with the existence of model completions/companions; see
\cite{Mac}.)  What about the dual question for the class of co-existentially
closed
members of {\bf CCH}?  (For example, if we could characterize the 
co-existentially closed
members of {\bf CCH} as the indecomposable curves, a co-elementary
class, then the answer to the last question would be yes, and we would have 
a nice example of a model complete co-elementary class of continua.)

$(viii)$ Two properties closely related to being co-existentially closed in the 
class of continua are class(W) and class(C).  A continuum $X$ is {\bf in
class(W)} (resp., {\bf in class(C)}) if whenever $f:Y \to X$ is a continuous
surjection and $Y$ is a continuum, then $f$ is weakly confluent (resp.,
confluent).  (These notions are due to A. Lelek \cite{Nad}.)  Clearly
if $X$ is a co-existentially closed continuum, then $X$ is in 
class(W).
Class(W) is definitely broader, since it contains all arcs. Spaces in
class(W) are {\bf unicoherent}; i.e., possessed of the feature that the
intersection of any two subcontinua whose union is everything must be
connected.  There are some nice characterizations of class(W), but they
go beyond our scope.  An elegant characterization of class(C) is simply
being hereditarily indecomposable.
 
$(ix)$ We would like to see a direct ultracoproduct proof of \ref{4.2},
one that doesn't rely on Banach model theory.


\begin{thebibliography}{99}
%\bibitem{Balb}
%     R. Balbes and P. Dwinger, ``Distributive Lattices,'' Univ. of 
%     Missouri Press, Columbia, 1974.
\bibitem{Bana}
     B. Banaschewski, ``More on compact Hausdorff spaces and finitary
     duality,'' {\sl Can. J. Math.\/} {\bf 36}(1984), 1113--1118.
%\bibitem{Ban1}
%     P. Bankston, ``Ultraproducts in topology,'' {\sl General Topology and its 
%     Applications\/} {\bf 7}(1977), 283--308.
\bibitem{Ban2}
     P. Bankston, ``Reduced coproducts of compact Hausdorff spaces,''
     {\sl J. Symbolic Logic\/} {\bf 52}(1987), 404--424. 
\bibitem{Ban3}
     -----------, ``Model-theoretic characterizations of arcs and simple
     closed curves,'' {\sl Proc. A. M. S.\/} {\bf 104}(1988), 898--904.
\bibitem{Ban4}
     ----------, ``Co-elementary equivalence for compact Hausdorff spaces
     and compact abelian groups,'' {\sl J. Pure and Applied Algebra\/}
     {\bf 68} (1990), 11--26.
\bibitem{Ban5}
     -----------, ``Taxonomies of model-theoretically defined topological
     properties,'' {\sl J. Symbolic Logic\/} {\bf 55}(1990), 589--603.
%\bibitem{Ban6}
%     -----------, ``Corrigendum to `Taxonomies of model-theoretically
%     defined topological properties','' {\sl J. Symbolic Logic\/}
%     {\bf 56}(1991), 425--426.
\bibitem{Ban8}
     -----------, ``Co-elementary equivalence, co-elementary maps, and
                    generalized arcs,'' {\sl Proc. A. M. S.} (to appear).
\bibitem{Ban7}
     -----------, ``On the topological model theory of normal disjunctive
                    lattices,'' (submitted).
%\bibitem{BS}
%     J. L. Bell and A. B. Slomson, ``Models and Ultraproducts,''
%     North Holland, Amsterdam, 1969.
\bibitem{CK}
     C. C. Chang and H. J. Keisler, ``Model Theory (Third Edition),'' 
     North Holland, Amsterdam, 1989.
%\bibitem{Coh}
%     P. M. Cohn, ``Universal Algebra,'' D. Reidel, Dordrecht, 1981.
%\bibitem{CN}
%     W. W. Comfort and S. Negrepontis, ``The Theory of Ultrafilters,''
%     Springer-Verlag, Berlin-New York, 1974.
\bibitem{DK}
     D. Dacunha-Castelle and J. Krivine, ``Applications des ultraproduits
     \`{a} l'\'{e}tude des espaces et des alg\`{e}bres de Banach,''
     {\sl Studia Math.\/} {\bf 41}(1972), 315--334.
\bibitem{Ekl}
     P. Eklof, ``Ultraproducts for algebraists,'' in ``Handbook of 
     Mathematical Logic,'' North Holland, Amsterdam, 1977, pp. 105--137.
\bibitem{E}
     R. Engelking, ``Outline of General Topology,'' North Holland, Amsterdam,
     1968.
%\bibitem{GJ}
%     L. Gillman and M. Jerison, ``Rings of Continuous Functions,'' Van
%     Nostrand, Princeton, 1960.
%\bibitem{Gra}
%     G. Gr\"{a}tzer, ``General Lattice Theory,'' Academic Press, New York,
%     1978.
\bibitem{Gur}
     R. Gurevi\v{c}, ``On ultracoproducts of compact Hausdorff spaces,''
     {\sl J. Symbolic Logic\/} {\bf 53}(1988), 294--300.
\bibitem{HH}
     S. Heinrich and C. W. Henson, ``Model theory of Banach spaces, II:
     isomorphic equivalence,'' {\sl Math. Nachrichten\/} {\bf 125}(1986),
     301--317.
\bibitem{HHM}
     S. Heinrich, C. W. Henson, and L. C. Moore, Jr., ``A note on elementary
     equivalence of $C(K)$ spaces,'' {\sl J. Symbolic Logic\/} {\bf 52}(1987),
     368--373.
\bibitem{Hen1}
     C. W. Henson, ``Nonstandard hulls of Banach spaces,'' {\sl Israel J.
     Math.\/} {\bf 25}(1976), 108--144.
\bibitem{Hen2}
     C. W. Henson, ``Nonstandard analysis and the theory of Banach spaces,''
     {\sl Springer Lecture Notes in Math., No. 983} (1983), 27--112.
\bibitem{HI} 
     C. W. Henson and J. Iovino, ``Banach Space Model Theory, I,''
     (lecture notes monograph, in preparation).
%\bibitem{HJRT}
%     C. W. Henson, C. G. Jockusch, Jr., L. A. Rubel, and G. Takeuti, 
%     ``First order topology,'' {\sl Dissertationes Mathematicae\/}
%     {\bf 143}(1977), 1--40.
%\bibitem{HY}
%     J. G. Hocking and G. S. Young, ``Topology,'' Addison-Wesley, Reading, MA,
%     1961. 
\bibitem{Hod}
      W. Hodges, ``Model Theory,'' Cambridge University Press, Cambridge,
      1993.
%\bibitem{Joh}
%     P. Johnstone, ``Stone Spaces,'' Cambridge University Press, Cambridge,
%     1982.
\bibitem{Kur}
     K. Kuratowski, ``Topology, vol. II,'' Academic Press, New York, 1968.
\bibitem{Mac}
     A. Macintyre, ``Model completeness,'' in ``Handbook of 
     Mathematical Logic,'' North Holland, Amsterdam, 1977, pp. 139--180.
%\bibitem{Mac}
%     S. Mac Lane, ``Categories for the Working Mathematician,'' Springer-
%     Verlag, Berlin-New York, 1971. 
%\bibitem{Mo}
%     R. L. Moore, ``Foundations of Point Set Theory,'' Amer. Math. Soc.
%     Colloq. Pub., vol. 13, Providence, RI, 1962.
\bibitem{Nad}
     S. B. Nadler, Jr., ``Continuum Theory, An Introduction,'' Marcel 
     Dekker, New York, 1992.
%\bibitem{Nag}
%     K. Nagami, ``Dimension Theory,'' Academic Press, New York, 1970.
%\bibitem{NTT}
%     J. Nikiel, H. M. Tuncali, and E. D. Tymchatyn, ``Continuous images of
%     arcs and inverse limit methods,'' {\sl Memoirs of the A. M. S.\/}
%     {\bf 104\/} (no. 498) (1993), 1--80.
%\bibitem{Ohk}
%     T. Ohkuma, ``Ultrapowers in categories,'' {\sl Yokohama Math. J.\/}
%     {\bf 14}(1966), 17--37.
%\bibitem{Pea}
%     A. R. Pears, ``Dimension Theory of General Spaces,'' Cambridge University
%     Press, Cambridge, 1975.
\bibitem{Q} W. V. O. Quine, ``Word and Object,'' The M. I. T. Press, 
      Cambridge, MA, 1960.
\bibitem{Ros}
     J. Rosick\'{y}, ``Categories of models,'' {\sl Seminarberichte Mathematik
     Informatik Fernuniversit\"{a}t\/} {\bf 19}(1984), 337--413.
\bibitem{She}
     S. Shelah, ``Every two elementarily equivalent models have
     isomorphic ultrapowers,'' {\sl Israel J. Math.} {\bf 10}(1971), 224--233.
\bibitem{Sim}
     G. F. Simmons, ``Introduction to Topology and Modern Analysis,''
     McGraw-Hill, New York, 1963.
%\bibitem{Wal}
%     R. C. Walker, ``The Stone-\v{C}ech Compactification,'' Springer-Verlag,
%     Berlin-New York, 1974.
\bibitem {Walm}
     H. Wallman, ``Lattices and topological spaces,'' {\sl Ann. Math.(2)\/}
     {\bf 39}(1938), 112--126.
\bibitem{Wil}
     S. Willard, ``General Topology,'' Addison-Wesley, Reading, MA, 1970.
\end{thebibliography}
\end{document}